\documentclass{article}
\usepackage{latexsym, color, amsmath,amsfonts,amssymb,graphicx,mathrsfs}
\usepackage[margin=1.3in]{geometry}
\numberwithin{equation}{section}
\newtheorem{theorem}{Theorem}[section]

\newtheorem{lemma}{Lemma}[section]
\newtheorem{proposition}{Proposition}[section]
\newtheorem{remark}{Remark}[section]
\newtheorem{example}{Example}[section]
\newtheorem{corollary}{Corollary}[section]
\def\proof{\mbox {\it Proof.~}}

\newcommand{\cqfd}{\mbox{}\nolinebreak\hfill\rule{2mm}{2mm}\medbreak\par}


\title{Non-potential systems with relativistic operators and maximal monotone boundary conditions}
\author{\\ Petru JEBELEAN$^{a,}$\footnote{petru.jebelean@e-uvt.ro - corresponding author}$\ $ and C\u{a}lin \c{S}ERBAN$^{b,}$\footnote{calin.serban@e-uvt.ro}\\ [0.10ex]\small
$^a$Institute for Advanced Environmental Research, $^b$Department of Mathematics\\ \small West University of Timi\c{s}oara,\\
\small Blvd. V. P\^{a}rvan, no. 4, 300223 - Timi\c{s}oara, Romania}
\date{}

\begin{document}
\maketitle

\begin{abstract}
\noindent We are concerned with solvability of a non-potential system involving two relativistic operators, subject to boundary conditions expressed in terms of maximal monotone operators. The approach makes use of a fixed point formulation and relies on a priori estimates and convergent to zero matrices.
\end{abstract}


\bigskip

\section{Introduction}
We deal with solvability of a system of type
\begin{equation}\label{eqpb1s}
	\left \{
            \begin{array}{ll}
		      -\left[ \phi_n(u^{\prime}) \right] ^{\prime}
              = f_1(t,u,v) \\
		      -\left[ \phi_m(v^{\prime}) \right] ^{\prime}
              = f_2(t,u,v)
			\end{array}  \quad \mbox{ in }[0,T],
\right.
\end{equation}
associated with the multivalued boundary conditions
\begin{equation}\label{bvpb1s}
\left \{
            \begin{array}{ll}
		      \left ( \phi_n \left( u^{\prime }\right)(0), -\phi_n \left( u^{\prime }\right)(T)\right )
              \in \gamma(u(0), u(T))\\
		      \left ( \phi_m \left( v^{\prime }\right)(0), -\phi_m \left( v^{\prime }\right)(T)\right )
              \in \eta(v(0), v(T))
			\end{array},
\right.
\end{equation}
where, for $q\in \{ n,m \}$, the homeomorphism $\phi_q:\mathbb{B}_1(q) \to \mathbb{R}^q$ is given by
\begin{equation*}\label{phi}
\displaystyle \phi_q(y)=\frac{y}{\sqrt{1- |y|^2}}\quad (y\in \mathbb{B}_1(q));
\end{equation*}
here and hereafter, $\mathbb{B}_r(q)$ denotes the open ball of radius $r$ in $\mathbb{R}^q$ endowed with the Euclidean norm $|\cdot|$ corresponding to the usual inner product $\langle \cdot | \cdot \rangle$. Also, on the product space $\mathbb{R}^q \times \mathbb{R}^q$ we consider the usual inner product $\langle \!  \langle  \cdot | \cdot \rangle \!  \rangle$.
The operators $\gamma: \mathbb{R}^n \times \mathbb{R}^n \to 2^{\mathbb{R}^n \times \mathbb{R}^n}$ and $\eta: \mathbb{R}^m \times \mathbb{R}^m \to 2^{\mathbb{R}^m \times \mathbb{R}^m}$ are maximal monotone, with $0_{\mathbb{R}^n \times \mathbb{R}^n}\in \gamma(0_{\mathbb{R}^n \times \mathbb{R}^n})$,  $0_{\mathbb{R}^m \times \mathbb{R}^m}\in \eta(0_{\mathbb{R}^m \times \mathbb{R}^m})$ and the mappings $f_{1}:[0,T]\times\mathbb{R}^n \times \mathbb{R}^m \to \mathbb{R}^n$, $f_{2}:[0,T]\times\mathbb{R}^n \times \mathbb{R}^m \to \mathbb{R}^m$ are continuous.
The notations $D(\gamma)$, $D(\eta)$ will stand for the domain of $\gamma$ and $\eta$, respectively (recall, this means $D(\gamma)= \{z\in \mathbb{R}^n \times \mathbb{R}^n \, : \, \gamma(z) \neq \emptyset \}$; similar definition for $D(\eta)$) and
we denote $C_q^1:=C^{1}(\left[ 0,T \right]; \mathbb{R}^q)$ ($q\in \{ n,m \}$).
\medskip

By \textit{a solution} of problem \eqref{eqpb1s}-\eqref{bvpb1s} we mean a couple of functions $(u,v)\in C_n^{1} \times C_m^1$ with  $| u^{\prime}(t)|<1$, $| v^{\prime}(t)|<1$ for all $t\in  [0,T]$, such that $\phi_n(u^{\prime })\in C_n^{1}$, $\phi_m(v^{\prime })\in C_m^{1}$, $(u(0), u(T))\in D\left(\gamma\right)$, $(v(0), v(T))\in D\left(\eta\right)$ and which satisfies \eqref{eqpb1s} and \eqref{bvpb1s}.
\medskip

In recent years a special attention was paid to various qualitative aspects for boundary value problems involving the {\it relativistic operator}: $u\mapsto\left[\phi_q(u')\right]'$, which is known that occurs in the dynamics of special relativity. Mainly, the results obtained for this operator use both topological techniques as well as variational methods. Thus, among others and far from being exhausted, related to existence and multiplicity of solutions of nonlinear systems or scalar equations subject to classical boundary conditions, such as Dirichlet, Neumann or periodic, under various assumptions on the perturbing nonlinearity, we refer the reader to e.g. \cite{[AmAr]} - \cite{[BrMa2]}, \cite{[JeMaSe-PAMS]}, \cite{[JePr2]}, \cite{[JeSe]} and the references therein.
This work is mainly motivated by paper \cite{[Je]}, which is concerned with a system having the form
\begin{equation}
-\left[ \phi_n(u^{\prime}) \right] ^{\prime}=\nabla_u F(t,u) \quad \mbox{in }[0,T],  \label{syspnew}
\end{equation}
subject to the boundary condition
\begin{equation}\label{bcsyspnew}
 \left ( \phi_n \left( u^{\prime }\right)(0), -\phi_n \left( u^{\prime }\right)(T)\right )\in \partial j(u(0), u(T)),
\end{equation}
under the assumptions that $j:\mathbb{R}^{n} \times \mathbb{R}^n \to (-\infty , +\infty]$ is a proper lower semicontinuous convex function with $j(0_{\mathbb{R}^n \times \mathbb{R}^n})=0$ and $0_{\mathbb{R}^n \times \mathbb{R}^n}\in \partial j (0_{\mathbb{R}^n \times \mathbb{R}^n})$; here, $\partial j$ stands for the subdifferential of $j$ in the sense of convex analysis \cite{[Rock]}. The mapping  $F =  F(t,u)  :   [0,T] \times \mathbb{R}^n \to \mathbb{R}$ is supposed to be continuous, with $\nabla_u F$ continuous on the set $[0,T] \times \mathbb{R}^n$ and $F(\, \cdot \, , 0_{\mathbb{R}^n})=0.$ Notice that by the lower semicontinuity assumption on the proper convex function $j$, the multivalued operator $\partial j : \mathbb{R}^n \times \mathbb{R}^n \rightarrow 2^{\mathbb{R}^n \times \mathbb{R}^n}$ is maximal monotone \cite{[Rock2]}. It is known (see e.g. \cite{[GaPa]} - \cite{[Je1]}, \cite{[KaPa]} and the references therein) that different choices of $j$ yield various boundary conditions, including the classical ones (Dirichlet, Neumann, periodic, antiperiodic, etc.) as well as other ones of special interest.
\smallskip

Observe that problem \eqref{syspnew}-\eqref{bcsyspnew} is entirely of potential type, both through the system and through the boundary condition. As shown in \cite{[Je]}, this allows a variational formulation in the frame of critical point theory for convex, lower semicontinuous perturbations of $C^1$-functionals. Then, taking the advantage of this key feature, therein one obtains the existence of minimum energy as well as saddle-point solutions of the problem. It is worth to notice that, as seen by concrete examples of applications, the generality of \eqref{syspnew}-\eqref{bcsyspnew} is quite broad. However, it is clear that firstly the potentiality of the system \eqref{syspnew} and secondly of the boundary condition \eqref{bcsyspnew} restrict the applicability of the results from \cite{[Je]}. In this view, it appears as being natural to consider a non-potential system
\begin{equation}
-\left[ \phi_n(u^{\prime}) \right] ^{\prime}=f(t,u) \quad \mbox{in }[0,T],  \label{syspnew1}
\end{equation}
associated with the non-potential boundary condition
\begin{equation}\label{bcsyspnew1}
 \left ( \phi_n \left( u^{\prime }\right)(0), -\phi_n \left( u^{\prime }\right)(T)\right )\in \gamma(u(0), u(T)),
\end{equation}
where the maximal monotone operator $\gamma$ is as above and $f:[0,T] \times \mathbb{R}^n \to \mathbb{R}^n$ is continuous. For the sake of clarity of the fact that \eqref{bcsyspnew1}  is really more general than \eqref{bcsyspnew}, we give here a simple example of a boundary condition which is expressed by a maximal monotone operator  that is not of subdifferential type. In this respect, let $A$ be a real ($2n \times 2n$)-matrix which is positive semi-definite (i.e. $\langle \! \langle A z \, | \, z  \rangle \! \rangle \geq 0$, for all $z\in \mathbb{R}^{n} \times \mathbb{R}^n$) and
\begin{equation}\label{gammaa}
\gamma_A (z) = A\, z \qquad (z\in \mathbb{R}^{n} \times \mathbb{R}^n).
\end{equation}
Then $\gamma_A$ is maximal monotone \cite[Example 1.5 (b)]{[Ph]}. If $A$ is not symmetric, then $\gamma_A$ will not be cyclically monotone (see p. 240 in \cite{[Rock]}) and thus, by a classical result on the potentiality of maximal monotone operators \cite{[Rock2]}, there is no lower semicontinuous proper convex function $j:\mathbb{R}^{n} \times \mathbb{R}^n \to (-\infty , +\infty]$ such that $\gamma_A = \partial j$. Summarizing, if $A$ is positive semi-definite, but it is not symmetric, the boundary condition
\begin{equation}\label{bcsysp}
\left(
\begin{array}{clrc}
  \phi_n \left( u^{\prime }\right)(0) \\[6pt]
  -\phi_n \left( u^{\prime }\right)(T)
\end{array}
\right)=A\
\left(
\begin{array}{clrc}
  u(0) \\[6pt]
  u(T)
\end{array}
\right)
\end{equation}
is of type \eqref{bcsyspnew1}, with $\gamma$ which is not the subdifferential of a proper lower semicontinuous convex function $j$.
\medskip

Coming back to system \eqref{eqpb1s}-\eqref{bvpb1s}, let us note that it contains as a particular case  problem \eqref{syspnew1}-\eqref{bcsyspnew1}. Indeed, given be $f$ and $\gamma$, this can be seen by taking $m=n$, $\eta = \gamma $, $f_1(t,u,v)=f(t,u)$ and $f_2(t,u,v)=f( t,v)$. Therefore, a qualitative result about \eqref{eqpb1s}-\eqref{bvpb1s} can be easily transposed to problem \eqref{syspnew1}-\eqref{bcsyspnew1} (see Corollaries \ref{CBS1} and \ref{apriori}). On the other hand, our main result (Theorem \ref{th2multma}) is of  type of \cite[Theorem 2.1]{[JePr1]} which is concerned with classical vector $p$ and $q$-Laplacian operators. In this view, the brought novelties consist in the fact that the differential operators in  \eqref{eqpb1s} are singular and the boundary conditions \eqref{bvpb1s} are of  non-potential type. As we will see, both technically as well as structurally these bring specific approaches and features, even if broadly speaking, it is about a similar strategy to the one in \cite{[JePr1]} - this makes use of a priori estimates combined with convergent to zero matrices. We mention that unlike the refered classical case, for system \eqref{eqpb1s}  subject to a whole class of boundary conditions of type \eqref{bvpb1s}, we obtain an "universal" existence result, meaning that no additional assumptions  must be imposed on the nonlinearities $f_1$, $f_2$  to ensure the existence of the solution (see Theorem \ref{th2multma} $(i)$).  Another notable aspect is that various choices of operators $\gamma$ and $\eta$ in \eqref{bvpb1s} allow the coupling of boundary conditions of different types.
\medskip

The rest of the paper is organized as follows. In Section 2 we give the fixed point formulation of problem \eqref{eqpb1s}-\eqref{bvpb1s}. Then, in Section 3 we present the main result concerning solvability of \eqref{eqpb1s}-\eqref{bvpb1s} and derive existence results for system \eqref{syspnew1}-\eqref{bcsyspnew1}.

\section{A fixed point approach}\label{sectiunea2}

Let $q\in\{n,m\}$. By $\| \cdot \| _{\infty}$ we denote the sup-norm on $C_q:=C([0,T]; \mathbb{R}^q)$ and the usual norm on $L^r:=L^r([0,T]; \mathbb{R}^q)$ will be denoted by $\| \cdot \|_{L^r}$ ($1\leq r \leq \infty$). The Sobolev space $H^1_q:=W^{1, 2}([0,T]; \mathbb{R}^q)$ will be considered with the norm
$$\|w\|_{H^1_q}=\left(\|w'\|_{L^2}^2+\|w\|_{L^2}^2\right)^{1/2}.$$
The product spaces $\mathcal{C}:=C_n\times C_m$ and $\mathcal{H}:=H^1_n\times H^1_m$ will be considered with the usual product topology; we choose to see $\mathcal{H}$ as being endowed with the norm
$$\|(u,v)\|_{\mathcal{H}}=\|u\|_{H^1_n}+\|v\|_{H^1_m} \quad ( (u,v)\in\mathcal{H}).$$

According to \cite[Theorem 2.8]{[Je]}, we have that for any $h\in C_n$, the system
\begin{equation}\label{eqpb2}
    -\left[ \phi_n(u^{\prime}) \right] ^{\prime} +u = h(t)\quad \mbox{in } [0,T],
\end{equation}
has an unique solution satisfying the boundary condition \eqref{bcsyspnew1}.
This enables us to define the solution operator
$S_{\gamma}:C_n\to H^1_n$ by
$S_{\gamma}(h):=\mbox{the unique solution of \eqref{eqpb2}-\eqref{bcsyspnew1}}\,  (h\in C_n).$
In a similar way  we introduce $S_{\eta}:C_m\to H^1_m$ by $S_\eta(l):=\mbox{the unique solution of}$ problem
\begin{equation*}\label{eqpb22}
    -\left[ \phi_m(v^{\prime}) \right] ^{\prime} +v = l(t)\ \ \mbox{in } [0,T],\quad \left ( \phi_m \left( v^{\prime }\right)(0), -\phi_m \left( v^{\prime }\right)(T)\right )\in \eta(v(0), v(T)) \qquad (l\in C_m).
\end{equation*}

\begin{proposition}\label{Scont}
The operators $S_{\gamma}$ and $S_\eta$ are continuous.
\end{proposition}
\proof We refer to $S_\gamma$; similar argument for $S_\eta$. Let $\{h_k\}\subset C_n$ be a sequence such that $h_k\to h\in C_n$ as $k \to \infty$, and set $u:=S_{\gamma}(h)$, $u_k:=S_{\gamma}(h_k)$ ($k\in \mathbb{N}$). As $\gamma$ is (maximal) monotone, by virtue of \eqref{bcsyspnew1} one has
\begin{eqnarray}\label{inegop}
  0 &\leq& \langle \!  \langle  \left ( \phi_n \left( u_k^{\prime }\right)(0),-\phi_n \left( u_k^{\prime }\right)(T)\right)-\left(\phi_n \left( u^{\prime }\right)(0),-\phi_n \left( u^{\prime }\right)(T)\right ) | \left(u_k(0),u_k(T)\right)- \left(u(0),u(T)\right)\rangle \!  \rangle \nonumber\\
   &=& \langle \!  \langle  \left ( \phi_n \left( u_k^{\prime }\right)(0)-\phi_n \left( u^{\prime }\right)(0),-\phi_n \left( u_k^{\prime }\right)(T)+\phi_n \left( u^{\prime }\right)(T)\right ) | \left(u_k(0)-u(0),u_k(T)-u(T)\right)\rangle \!  \rangle \nonumber\\
   &=& \langle\phi_n \left( u_k^{\prime }\right)(0)-\phi_n \left( u^{\prime }\right)(0) | u_k(0)-u(0)\rangle - \langle\phi_n \left( u_k^{\prime }\right)(T)-\phi_n \left( u^{\prime }\right)(T) | u_k(T)-u(T)\rangle.
\end{eqnarray}

On the other hand, multiplying the equality
$$-\left[\phi_n(u_k^{\prime})- \phi_n(u^{\prime})\right] ^{\prime}+u_k-u=h_k(t)-h(t) $$
by $u_k-u$, integrating over $[0,T]$ and using the integration by parts formula, we derive
$$\langle\phi_n \left( u_k^{\prime }\right)(0)-\phi_n \left( u^{\prime }\right)(0) | u_k(0)-u(0)\rangle - \langle\phi_n \left( u_k^{\prime }\right)(T)-\phi_n \left( u^{\prime }\right)(T) | u_k(T)-u(T)\rangle+$$
$$+\int_0^T\langle\phi_n \left( u_k^{\prime }\right)-\phi_n \left( u^{\prime }\right) | u_k^{\prime}-u'\rangle+\|u_k-u\|_{L^2}^2=\int_0^T\langle h_k(t)-h(t) | u_k-u\rangle$$
and from \eqref{inegop}, one obtains
$$\int_0^T\langle\phi_n \left( u_k^{\prime }\right)-\phi_n \left( u^{\prime }\right) | u_k^{\prime}-u'\rangle+\|u_k-u\|_{L^2}^2\leq\sqrt{T}\|h_k-h\|_\infty\|u_k-u\|_{L^2}.$$
Then, since $\phi_n$ satisfies \cite[Lemma 2.1]{[JeMaSe-PAMS]}:
$$\langle\phi_n(x)-\phi_n(y) | x-y\rangle\geq|x-y|^2\quad (x,y\in \mathbb{B}_1(n)),$$
we get that $\{\|u_k-u\|_{L^2}\}$ is bounded and
$$\|u_k-u\|_{H^1_n}^2\leq\sqrt{T}\|h_k-h\|_\infty\|u_k-u\|_{L^2},$$
which, as $h_k\to h$ in $C_n$, yields that $u_k\to u$ in $H^1_n$ and the continuity of $S_{\gamma}$ is proved.\cqfd
\medskip

Next, we introduce the mappings $g_1:[0,T]\times\mathbb{R}^n \times \mathbb{R}^m \to \mathbb{R}^n$, $g_2:[0,T]\times\mathbb{R}^n \times \mathbb{R}^m \to \mathbb{R}^m$ by
$$g_1(t,x,y)=f_1(t,x,y)+x,\ \quad g_2(t,x,y)=f_2(t,x,y)+y \quad \, (t\in[0,T], \, x\in \mathbb{R}^n,\ y\in\mathbb{R}^m)$$
and define the \textit{Nemytskii operators} $N_{g_1}:\mathcal{C}\to C_n$ and $N_{g_2}:\mathcal{C}\to C_m$ by
$$N_{g_1}(u,v)(t)=g_1(t,u(t),v(t)),\ \quad  N_{g_2}(u,v)(t)=g_2(t,u(t),v(t))\quad \, (t\in[0,T],\ (u,v)\in \mathcal{C}).$$
It is well known that $N_{g_1}$, $N_{g_2}$ are continuous and bounded (i.e., takes bounded sets from $\mathcal{C}$ into bounded sets in $C_q$ ($q\in\{n,m\}$)). Then, let $Q_\gamma:\mathcal{H}\to H^1_n$ and $Q_\eta:\mathcal{H}\to H^1_m$ be given by
\begin{equation}\label{qdef}
    Q_\gamma=S_{\gamma}\circ N_{g_1}\circ i_d, \quad Q_\eta=S_{\eta}\circ N_{g_2}\circ i_d,
\end{equation}
where $i_d:\mathcal{H}\to \mathcal{C}$ is the identity map.

\begin{proposition}\label{Qcompfix}
$(i)$ The operator $Q=(Q_\gamma,Q_\eta):\mathcal{H}\to\mathcal{H}$ is compact.

\smallskip
\noindent $(ii)$ A couple of functions $(u,v)\in \mathcal{H}$ is a solution of problem \eqref{eqpb1s}-\eqref{bvpb1s} iff it is a fixed point of $Q$.
\end{proposition}
\proof Statement $(i)$ follows from Proposition \ref{Scont}, the continuity of $N_{g_1}$, $N_{g_2}$ and the compactness of $i_d$, while $(ii)$ is straightforward.\cqfd
\medskip

Before concluding this section, we note for reader convenience some simple characterizations of the convergent to zero matrices, which will be needed in the sequel. Recall, a square matrix $M$ is said to be \textit{convergent to zero} if $M^k\to 0$ as $k\to\infty$. The following result is proved in \cite[Lemma 2]{[Pr]}.

\begin{lemma}\label{matrixconv0}
Let $M$ be a square matrix of non-negative numbers. The following
statements are equivalent:
\begin{enumerate}
  \item[$($i$)$] $M$ is a convergent to zero matrix;
  \item[$($ii$)$] $I-M$ is non-singular and
               $$(I-M)^{-1}=I+M+M^2+\ldots;$$
  \item[$($iii$)$] $|\lambda|<1$ for every $\lambda\in\mathbb{C}$ with $\det(M-\lambda I)=0$;
  \item[$($iv$)$] $I-M$ is non-singular and $(I-M)^{-1}$ has non-negative elements.
\end{enumerate}
\end{lemma}

\section{Main result}\label{sectiunea3}
We denote
\begin{equation*}\label{defK}
    \mathcal{K}(q):= \left \{ w \in W^{1, \infty}([0,T]; \mathbb{R}^q) \, : \, \|w^{\prime}\|_{_{L^{\infty}}}\leq 1 \right \}\quad (q\in \{n,m \})
\end{equation*}
and define the first eigenvalue-like constants:
\begin{equation}\label{lambda1n}
   \lambda_1(\gamma):=\inf\left \{ \frac{\|u^{\prime} \|_{L^2}^2}{\|u \|_{L^2}^2} \, : \, u\in \mathcal{K}(n)\setminus \{ 0_{\mathbb{R}^n} \}, \, (u(0),u(T))\in  D(\gamma) \right \},
\end{equation}
\begin{equation}\label{lambda1m}
   \lambda_1(\eta):=\inf\left \{ \frac{\|v^{\prime} \|_{L^2}^2}{\|v \|_{L^2}^2} \, : \, v\in \mathcal{K}(m)\setminus \{ 0_{\mathbb{R}^m} \}, \, (v(0),v(T))\in  D(\eta) \right \}.
\end{equation}

The following main result provides sufficient conditions which ensures the solvability of the boundary value problem \eqref{eqpb1s}-\eqref{bvpb1s}.

\begin{theorem}\label{th2multma} Assume that one of the following four conditions is satisfied:

\medskip
\noindent $(i)$ $\lambda_1(\gamma)>0$, $\lambda_1(\eta)>0;$

\smallskip
\noindent $(ii)$ $\lambda_1(\eta)>0$ and there are constants $a\in [0,1), b, \delta \in \mathbb{R}_+$ such that
\begin{equation}\label{inegalF1ma}
    \langle f_1(t,x,y) \, | \, x  \rangle \leq (a-1)|x|^2+b|y|^2+\delta \quad (t\in[0,T],\ x \in \mathbb{R}^n,\ y \in \mathbb{R}^m);
\end{equation}

\smallskip
\noindent $(iii)$ $\lambda_1(\gamma)>0$ and there are constants $c, \delta\in \mathbb{R}_+, d\in [0,1)$  such that
\begin{equation}\label{inegalF1maa}
    \langle f_2(t,x,y) \, | \, y  \rangle \leq c|x|^2+(d-1)|y|^2+\delta \quad (t\in[0,T],\ x \in \mathbb{R}^n,\ y \in \mathbb{R}^m);
\end{equation}

\smallskip
\noindent $(iv)$ there are constants $a,b,c,d, \delta\in \mathbb{R}_+$ such that the matrix
\begin{equation}  \label{systmmug}
M:= \left[
\begin{array}{clrc}
a & b \\
c & d
\end{array}\right]
\end{equation}
is convergent to zero and
\begin{equation}  \label{syst_mug}
	\begin{array}{ll}
		\langle f_1(t,x,y) \, | \, x  \rangle
		\leq (a-1)|x|^2+b|y|^2+\delta \\
		\langle f_2(t,x,y) \, | \, y  \rangle
		\leq c|x|^2+(d-1)|y|^2+\delta
	\end{array} \quad (t\in[0,T],\ x \in \mathbb{R}^n,\ y \in \mathbb{R}^m).
\end{equation}

\noindent Then problem \eqref{eqpb1s}-\eqref{bvpb1s} has at least one solution.
\end{theorem}
\proof On account of Proposition \ref{Qcompfix} $(ii)$, it suffices to show that under any of the assumptions $(i)-(iv)$, the operator $Q$ has a fixed point. In this view, according to Schaefer's theorem (see e.g. \cite[Corollary 4.4.12]{[L]}), it suffices to prove that the set
$$\{(u,v)\in \mathcal{H}:\ \exists\ \lambda\in(0,1]\mbox{ such that } (u,v)=\lambda Q(u,v)\}$$
is bounded in $(\mathcal{H},\|\cdot\|_{\mathcal{H}})$.
\smallskip

Let $(u,v)\in \mathcal{H}$ be such that $Q(u,v)=\lambda^{-1}(u,v)$, with some $\lambda\in(0,1]$. From \eqref{qdef} and definition of operators $S_{\gamma}$ and $S_{\eta}$, this means that $|\lambda ^{-1}u'(t)|<1$, $|\lambda ^{-1}v'(t)|<1$ for all $t\in[0,T]$, $(\lambda^{-1}u(0), \lambda^{-1}u(T))\in D\left(\gamma\right)$, $(\lambda^{-1}v(0), \lambda^{-1}v(T))\in D\left(\eta\right)$,
\begin{equation}\label{apr1n}
-\left[ \phi_n\left(\lambda^{-1}u^{\prime}\right) \right] ^{\prime} +\lambda^{-1}u = g_1(t,u,v)\quad \, (t\in[0,T]),
\end{equation}
\begin{equation}\label{apr1m}
-\left[ \phi_m\left(\lambda^{-1}v^{\prime}\right) \right] ^{\prime} +\lambda^{-1}v = g_2(t,u,v)\quad \, (t\in[0,T])
\end{equation}
and
\begin{equation}\label{apr2n}
    \left ( \phi_n \left(\lambda^{-1}u^{\prime }\right)(0), -\phi_n \left(\lambda^{-1}u^{\prime }\right)(T)\right )\in \gamma(\lambda^{-1}u(0), \lambda^{-1}u(T)),
\end{equation}
\begin{equation}\label{apr2m}
    \left ( \phi_m \left(\lambda^{-1}v^{\prime }\right)(0), -\phi_m \left(\lambda^{-1}v^{\prime }\right)(T)\right )\in \eta(\lambda^{-1}v(0), \lambda^{-1}v(T)).
\end{equation}

\noindent$(i)$ From \eqref{lambda1n} and \eqref{lambda1m}, we have that
$$\|u\|_{L^2}^2\leq\frac{\|u'\|_{L^2}^2}{\lambda_1(\gamma)}<\frac{T}{\lambda_1(\gamma)}\ \mbox{ and }\ \|v\|_{L^2}^2\leq\frac{\|v'\|_{L^2}^2}{\lambda_1(\eta)}<\frac{T}{\lambda_1(\eta)},$$
which gives $\|(u,v)\|_{\mathcal{H}}<\left(T+T/\lambda_1(\gamma)\right)^{1/2}+\left(T+T/\lambda_1(\eta)\right)^{1/2}$.
\medskip

\noindent $(ii)$ Multiplying \eqref{apr1n} by $\lambda^{-1}u$, then integrating over $[0,T]$ and using integration by parts formula, we infer
\begin{equation*}
    \left\langle\phi_n \left(\lambda^{-1}u^{\prime}\right)(0) | \lambda^{-1}u(0)\right\rangle - \left\langle\phi_n \left(\lambda^{-1}u^{\prime}\right)(T) | \lambda^{-1}u(T)\right\rangle+
\end{equation*}
\begin{equation}\label{apr3n}
    +\int_0^T\left\langle\phi_n \left(\lambda^{-1}u^{\prime}\right) | \lambda^{-1}u'\right\rangle+\lambda^{-2}\int_0^T|u|^2=\lambda^{-1}\int_0^T\langle g_1(t,u,v) | u\rangle.
\end{equation}
Using that $0_{\mathbb{R}^n \times \mathbb{R}^n}\in \gamma(0_{\mathbb{R}^n \times \mathbb{R}^n})$, the monotonicity of $\gamma$ and \eqref{apr2n}, one has
$$\left\langle\phi_n \left(\lambda^{-1}u^{\prime}\right)(0) | \lambda^{-1}u(0)\right\rangle - \left\langle\phi_n \left(\lambda^{-1}u^{\prime}\right)(T) | \lambda^{-1}u(T)\right\rangle\geq0,$$
which, on account of $\langle \phi_n(y)|y\rangle \geq |y|^2$ ($y\in \mathbb{B}_1(n)$) and \eqref{apr3n}, yields
\begin{equation}\label{apr4n}
    \|u\|_{H^1_n}^2\leq\lambda\int_0^T\langle g_1(t,u,v) | u\rangle = \lambda \int_0^T\langle f_1(t,u,v) | u\rangle+\lambda\|u\|_{L^2}^2.
\end{equation}
Then, from \eqref{lambda1m}, since $\lambda_1(\eta)>0$, one has $\|v\|_{L^2}^2< T/\lambda_1(\eta)$ and by virtue of \eqref{inegalF1ma} and \eqref{apr4n}, we get
$$\|u\|_{H^1_n}^2\leq a\|u\|_{L^2}^2+b\|v\|_{L^2}^2+\delta T<a\|u\|_{H^1_n}^2+k_1,$$
with $k_1=k_1(\eta):=(b/\lambda_1(\eta)+\delta)T$. Hence,
%$\|u\|_{L^2}^2< k_1/(1-a)$ and
$\|(u,v)\|_{\mathcal{H}}<\left(k_1/(1-a)\right)^{1/2}+\left(T+T/\lambda_1(\eta)\right)^{1/2}$.
\medskip

\noindent $(iii)$ Similarly to $(ii)$, using \eqref{apr1m}, \eqref{apr2m} and the monotonicity of $\eta$, we obtain
\begin{equation}\label{apr4m}
    \|v\|_{H^1_m}^2\leq\lambda\int_0^T\langle g_2(t,u,v) | v\rangle\leq \lambda \int_0^T\langle f_2(t,u,v) | v\rangle+\lambda\|v\|_{L^2}^2,
\end{equation}
which, together with \eqref{lambda1n} and \eqref{inegalF1maa}, imply $\|(u,v)\|_{\mathcal{H}}<\left(T+T/\lambda_1(\gamma)\right)^{1/2}+\left(k_2/(1-d)\right)^{1/2},$
where $k_2=k_2(\gamma):=(c/\lambda_1(\gamma)+\delta)T$.
\medskip

\noindent $(iv)$ From \eqref{syst_mug}, \eqref{apr4n} and \eqref{apr4m}, we derive
$$\|u\|_{H^1_n}^2\leq a\|u\|_{H^1_n}^2+b\|v\|_{H^1_m}^2+\delta T\ \mbox{ and }\
\|v\|_{H^1_m}^2\leq c\|u\|_{H^1_n}^2+d\|v\|_{H^1_m}^2+\delta T,$$
which can be written in the following vector form:
\begin{equation}\label{matricial}
    (I-M)\left[
    \begin{array}{c}
    \|u\|_{H^1_n}^2 \\[6pt]
    \|v\|_{H^1_m}^2
    \end{array}
    \right]\leq
    \left[
    \begin{array}{c}
    \delta T \\[6pt]
    \delta T
    \end{array}
    \right].
\end{equation}
Then, since matrix $M$ is convergent to zero, Lemma \ref{matrixconv0} guarantees that $I-M$ is invertible and the elements of $(I-M)^{-1}$ are all non-negative. Therefore, we may multiply \eqref{matricial} by $(I-M)^{-1}$ without changing the sense of the inequality and so, we get
\begin{equation*}
   \left[
    \begin{array}{c}
    \|u\|_{H^1_n}^2 \\[6pt]
    \|v\|_{H^1_m}^2
    \end{array}
    \right]
    \leq (I-M)^{-1}
        \left[
    \begin{array}{c}
    \delta T \\[6pt]
    \delta T
    \end{array}
    \right],
\end{equation*}
which shows that $\|(u,v)\|_{\mathcal{H}}$ is bounded by a constant depending on  $a,b,c,d, \delta$ and independent of $\lambda\in(0,1]$.
\cqfd
\medskip

\begin{remark}\label{rempoz}{\em Denote by $cone \,D(\gamma)$ the conical hull of the set $D(\gamma)$, that is
$$cone \,D(\gamma) := \{ \alpha z \, : \, \alpha \geq 0 , \, z\in D( \gamma )\}.$$
A sufficient condition to have $\lambda_1(\gamma)>0$ is
\begin{equation}\label{lugpoz}
\overline{cone \,D(\gamma)} \cap \{ (\xi, \xi) \, : \, \xi\in \mathbb{R}^n \} = \{ 0_{\mathbb{R}^n \times \mathbb{R}^n} \}.
\end{equation}
To see this, we argue as in the proof of Corollary 4.2 from \cite{[Je]}. Thus, setting
${\mathcal M}(n):= \{u\in H^1_n \, : \, ( u(0), u(T) ) \in \overline{cone \,D(\gamma)} \}$ and
$$\underline{\lambda}_1(\gamma ):= \inf \left \{ \frac{\|u^{\prime} \|_{L^2}^2}{\|u \|_{L^2}^2} \, : \, u\in {\mathcal M}(n) \setminus \{ 0_{\mathbb{R}^n} \} \right \},$$
from \cite[Theorem 3.1]{JePr11}, we know that $\underline{\lambda}_1(\gamma )>0$ and the conclusion follows from $\lambda_1(\gamma)\geq \underline{\lambda}_1(\gamma ).$
\smallskip

Similarly to \eqref{lugpoz}, condition
\begin{equation*}\label{eugpoz}
\overline{cone \,D(\eta)} \cap \{ (\xi, \xi) \, : \, \xi \in \mathbb{R}^m \} = \{ 0_{\mathbb{R}^m \times \mathbb{R}^m} \}
\end{equation*}
is a sufficient one to have that $\lambda_1(\eta)>0$.
}
\end{remark}
\smallskip

\begin{example}\label{primex}{\em Let $a,b \in \mathbb{R}$,  $h \in C_n$, $l \in C_m$ and  consider the system
\begin{equation}\label{expb1ss}
\left \{
            \begin{array}{ll}
		      -\displaystyle \left[ \phi_n(u') \right] ^{\prime}
              = (a-1)u + \displaystyle\frac{b}{1+|u|^2}\left [|v|^2 u +h(t) \right ] \\
		      -\displaystyle \left[ \phi_m(v') \right] ^{\prime}
              = \displaystyle\frac{b}{1+|v|^2}\left [|u|^2 v +l(t) \right ]+(a-1)v
			\end{array}  \quad \mbox{ in }[0,T].
\right.
\end{equation}
Then \eqref{expb1ss} associated with the Dirichlet and antiperiodic boundary conditions
\begin{equation*}\label{expb1s}
\left \{
            \begin{array}{ll}
		      u(0)=0_{\mathbb{R}^n}=u(T)\\
		      v(0)+v(T)=0_{\mathbb{R}^m}=v'(0)+v'(T)
			\end{array}
\right.
\end{equation*}
has at least one solution for any constants $a,b\in \mathbb{R}.$  This easily follows from Theorem \ref{th2multma} $(i)$ and Remark \ref{rempoz}, with

\begin{equation*} \label{dirant}
\gamma(x,y)=	\left \{
            \begin{array}{ll}
		      \mathbb{R}^n \times \mathbb{R}^n & \mbox{ if } x=0_{\mathbb{R}^n}=y  \\
		      \emptyset & \mbox{ otherwise, }
			\end{array}
\right. \qquad
\eta(x,y)=	\left \{
            \begin{array}{ll}
		      \{ (\xi, \xi) \, : \, \xi\in \mathbb{R}^m \} & \mbox{ if } x+y=0_{\mathbb{R}^m}  \\
		      \emptyset & \mbox{ otherwise. }
			\end{array}
\right.
\end{equation*}
On the other hand, if we take $\gamma(x,y)=(y,-x)$ ($(x,y) \in \mathbb{R}^n \times \mathbb{R}^n$), which is a maximal monotone operator \cite[Example 2.23 (a)]{[Ph]} of type $\gamma_A$ (see \eqref{gammaa}) with
$$A= \left[
\begin{array}{clrc}
\mathbb{O}_n & \mathbb{I}_n \\
-\mathbb{I}_n & \mathbb{O}_n
\end{array}\right] ,$$
 and
\begin{equation*} \label{dirant2}
\eta(x,y)=	\left \{
            \begin{array}{ll}
		      \{ (\xi, -\xi) \, : \, \xi\in \mathbb{R}^m \} & \mbox{ if } x=y  \\
		      \emptyset & \mbox{ otherwise, }
			\end{array}
\right.
\end{equation*}
then the boundary conditions \eqref{bvpb1s} become (see \eqref{bcsysp}):
\begin{equation}\label{expb1sss}
\left \{
            \begin{array}{ll}
		      \phi _n(u'(0))=u(T), \ \phi _n(u'(T))=u(0)\\
		      v(0)-v(T)=0_{\mathbb{R}^m}=v'(0)-v'(T).
			\end{array}
\right.
\end{equation}
In this case $\lambda_1(\gamma)=0=\lambda_1(\eta)$ and by Theorem \ref{th2multma} $(iv)$ we get that problem \eqref{expb1ss}-\eqref{expb1sss} has at least one solution provided that $a,b\in\mathbb{R}_+$ and $a+b<1$.
}
\end{example}

\begin{corollary}\label{CBS1}
Problem \eqref{syspnew1}-\eqref{bcsyspnew1} has at least one solution if one of the following two conditions is satisfied:

\medskip
\noindent $(i)$ $\lambda_1(\gamma)>0;$

\smallskip
\noindent $(ii)$ there are constants $a\in [0,1)$, $\delta \in \mathbb{R}_+$ such that
\begin{equation}\label{inegalF1mac}
    \langle f(t,x) \, | \, x  \rangle \leq (a-1)|x|^2+\delta, \quad \mbox{ for all }t\in[0,T], \, x \in \mathbb{R}^n.
\end{equation}
\end{corollary}
\proof Theorem \ref{th2multma} applies with $m=n$, $\eta=\gamma$, $f_1(t,u,v)=f(t,u)$ and $f_2(t,u,v)=f(t,v)$. If $(u,v)$ solves \eqref{eqpb1s}-\eqref{bvpb1s} with these choices of the data, then each of the functions $u$ or $v$ is a solution of \eqref{syspnew1}-\eqref{bcsyspnew1}. Clearly, \eqref{inegalF1mac} implies that \eqref{syst_mug} hold true with $b=c=0$ and $d=a$. The matrix $M$ in \eqref{systmmug} becomes
\begin{equation*}
M= \left[
\begin{array}{clrc}
a & 0 \\
0 & a
\end{array}\right]
\end{equation*}
and it is convergent to zero because $a\in[0,1).$ \cqfd

\begin{corollary}\label{apriori} If
\begin{equation}\label{fapriori}
 \limsup_{|x|\to\infty}\frac{\langle f(t,x) | x\rangle}{|x|^2}<0,\quad \mbox{uniformly with } t\in[0,T],
\end{equation}
then condition $(ii)$ in Corollary \ref{CBS1} is fulfilled, and hence problem \eqref{syspnew1}-\eqref{bcsyspnew1} has at least one solution.
\end{corollary}
\proof  From \eqref{fapriori} we can find constants $\sigma\in (0,1]$ and $\rho>0$ such that
\begin{equation*}
    \langle f(t,x) | x\rangle\leq-\sigma|x|^2,\quad \mbox{for all }\ t\in[0,T] \mbox{ and } x\in\mathbb{R}^n \mbox{ with } |x|>\rho.
\end{equation*}
Putting $k=k(\rho):= \max \{ | \langle f(t,x) | x\rangle |\, : \, (t,x)\in [0,T] \times \mathbb{B}_{\rho}(n) \},$ one has
\begin{equation*}
    \langle f(t,x) | x\rangle\leq-\sigma|x|^2+\sigma\rho^2+k,\quad \mbox{for all }\ t\in[0,T] \mbox{ and } x\in\mathbb{R}^n.
\end{equation*}
Then, \eqref{inegalF1mac} holds true with $a=1-\sigma$ and $\delta=\sigma\rho^2+k$. \cqfd

\end{document}